\documentclass{amsart}
\usepackage[mathscr]{eucal}
\usepackage{graphicx}
\usepackage{amscd}
\usepackage{amsfonts}
\usepackage{amsmath}
\usepackage{amsthm}
\usepackage{amssymb}
\usepackage{latexsym}

\newtheorem{theorem}[]{Theorem}
\newtheorem{proposition}[]{Proposition}

\newtheorem{lemma}[]{Lemma}
\theoremstyle{definition}                               

\newtheorem{remark}[]{Remark}

\def\question{\noindent\textbf{Question.} }
\def\remark{\noindent\textbf{Remark.} }

\newcommand{\beq}{\begin{equation}\label}

\newcommand{\into}{\,\,\hookrightarrow\,\,}

\newcommand{\onto}{\,\,\twoheadrightarrow\,\,}

\newcommand{\Spec}{{\mathtt{Spec}}}
\newcommand{\Hilb}{{\mathtt{Hilb}}}
\newcommand{\End}{{\mathtt{End}}}
\newcommand{\Hom}{{\mathtt{Hom}}}
\newcommand{\Aut}{{\mathtt{Aut}}}
\newcommand{\Mod}{{\mathtt{Mod}}}
\newcommand{\Ker}{{\mathtt{Ker}}}
\newcommand{\Ann}{{\mathtt{Ann}}}
\renewcommand{\dim}{{\mathtt{dim}}}
\renewcommand{\deg}{{\mathtt{deg}}}
\newcommand{\Par}{{\mathtt{Part}}}

\newcommand{\rk}{{\mathtt{rk}}}
\newcommand{\ch}{\chi}

\newcommand{\grd}{{\mathtt{gr}}}

\newcommand{\tr}{{{\mathtt{Tr}}}}

\newcommand{\dd}{{\mathcal{D}}}
\newcommand{\GL}{\mathtt{GL}}

\newcommand{\SSp}{\mathtt{Sp}}

\newcommand{\e}{{\mathbf{e}}}
\newcommand{\cc}{{\mathbf{c}}}

\newcommand{\B}{{\mathcal{H}}}

\newcommand{\triv}{{\mathtt {triv}}}

\newcommand{\sch}{\boldsymbol{\mathsf{s}}}
\newcommand{\mon}{\boldsymbol{\mathsf{m}}}

\newcommand{\h}{{\mathfrak{h}}}

\def\C{{\mathbb{C}}}
\def\Q{{\mathbb{Q}}}

\def\ms#1{\mathcal{#1}}

\def\Z{{\mathbb{Z}}}

\def\K{{\mathcal K}}

\def\sll2{{\mathfrak{s}\mathfrak{l}}_2}

\def\Cr{\C^{\mathsf{reg}}}

\def\H{{\mathbf{H}}}
\def\B{{\mathbf{B}}}

\def\HH{{\mathtt{HH}}}
\def\K{{\mathtt{K}}}

\def\ccirc{{{}_{^{^\circ}}}}

\begin{document}
\title{Morita Equivalence of Cherednik Algebras}
\author{Yuri Berest}
\address{Department of Mathematics,
Cornell University, Ithaca, NY 14853-4201, USA}
\email{berest@math.cornell.edu}
\author{Pavel Etingof}
\address{Department of Mathematics, MIT, Cambridge, MA 02139, USA}
\email{etingof@math.mit.edu}
\author{Victor Ginzburg}
\address{Department of Mathematics, University of Chicago,
Chicago, IL 60637, USA}
\email{ginzburg@math.uchicago.edu}
\begin{abstract}
We classify the rational Cherednik algebras $\, H_c(W) \,$ 
(and their spherical subalgebras) up to isomorphism and 
Morita equivalence in case when $\, W \,$ is the symmetric group
and $ c $ is a generic parameter value. 
\end{abstract}
\maketitle
\section{Introduction} 
In this paper, which is a sequel to our earlier work \cite{BEG},
we deal with Morita classification of rational Cherednik algebras. 
Our goal is to prove one of the main conjectures of \cite{BEG} 
(see Conjecture~8.12, {\it loc.~cit.}) in the case of symmetric groups.
Before stating our results we recall the notation and some of the
basic definitions.

Let $ W $ be a finite Coxeter group generated by reflections 
in a finite-dimensional complex vector space $ \h \,$. 
The {\it rational Cherednik algebras} associated to $ W $
is a family of associative algebras $ \{ H_c(W)\} $ 
parametrized by the set of $W$-invariant complex multiplicities
$\, c: R \to \C \,$ on the system of roots $ R \subset \h^* $ of 
$ W\,$.  Specifically, for a fixed $\, c: \alpha \mapsto 
c_{\alpha} \,$, the algebra $ H_c = H_c(W) $ is generated by the 
vectors of $\h \,$, $\,\h^*\,$, and the elements of $ W $
subject to the following relations  
$$
\begin{array}{lll}\displaystyle
&{}_{_{\vphantom{x}}}w\, x\, w^{-1}= w(x)\;\;,\;\;
w\, y\, w^{-1}= w(y)\,,&
\forall y\in \h\,,\,x\in \h^*\,,\,w \in W\, ,\break\medskip\\
&{}^{^{\vphantom{x}}}{}_{_{\vphantom{x}}}x_1\, x_2 = 
 x_2\, x_1\enspace,\enspace
y_1\, y_2=y_2\, y_1\,, &
\forall y_1,y_2\in \h,\;x_1\,,\,x_2 \in \h^*\, ,\break\medskip\\
&{}^{^{\vphantom{x}}}y\, x-x\, y = \langle y,x\rangle 
- \!\!\sum\limits_{\alpha\in R_{+}}\!\! c_\alpha\,\langle y,\alpha\rangle
\langle\alpha^\vee,x\rangle \, s_\alpha\,,& \forall y\in
\h\,,\,x\in \h^*\, .
\end{array}
$$
Here, as usual, we write $\, \alpha^\vee \in \h \,$ for the coroot, 
$\, s_\alpha \in W \,$ for the reflection corresponding to the 
root $\, \alpha \in R \,$, $\, R_{+} \subset R \,$ for a choice
of `positive' roots in $ R\,$, and 
$\,
\langle\,\boldsymbol{\cdot}\,,\,\boldsymbol{\cdot}\,\rangle \,$ for the canonical pairing
between $ \h $ and $ \h^* \,$. 

Note that the group algebra $ \C W $ of $ W $ embeds naturally
in $ H_c(W) \,$, and $\, H_0(W) = \dd(\h) \# W \,$, 
where $\,\dd(\h) \# W \,$ is the crossed product of $ \C W $ with
the ring of polynomial differential operators on $ \h \,$.
Furthermore, each $ H_c $ contains a distinguished subalgebra
$\, B_c := \e H_c \e \,$, where 
$\, \e := \frac{1}{|W|} \sum\,w\,$ is the 
symmetrizing idempotent in $\C W \subset H_c $ which
plays a role of the identity element in $ B_c \,$.
We call $ B_c = B_c(W) $ the {\it spherical algebra} associated 
to $\,(W,c)\,$. Since $\, B_0 = \e(\dd(\h) \# W)\e \cong \dd(\h)^W \,$,
the family $\, \{B_c(W)\} \,$ should be thought of as a deformation 
(in fact, the {\it universal} deformation) 
of the ring $ \dd(\h)^W $ of $W$-invariant differential 
operators on $\, \h\,$.

The structure of and the relationship between the algebras $ H_c $ and 
$ B_c $ depend drastically on the values of the parameter $ c\,$. It 
turns out that both of these are governed by the properties of the 
{\it Hecke algebra} $ \mathcal{H}_{W}(q) $ of $ W $ with 
$ q = e^{2\pi i c}\,$. 
In fact, one of the main results of \cite{BEG} states
\begin{theorem}[\cite{BEG}, Theorem~3.1, Corollary~4.2]
\label{T0}
If $ \mathcal{H}_{W}(e^{2\pi i c}) \,$ is semisimple, then
$ H_c $ and $ B_c $ are simple algebras,
Morita equivalent to each other.
\end{theorem}
Following \cite{BEG}, we will call $ c $ {\it regular} 
if the assumption of Theorem~\ref{T0} holds (i.e., if the 
Hecke algebra $ \mathcal{H}_{W}(e^{2\pi i c}) \,$ is semisimple).

In the present paper we will classify the algebras
$  H_c(W) $ and $ B_c(W) $ up to isomorphism
and Morita equivalence in case when $ W $ is a symmetric
group. Thus, from now on, we assume that $ W := S_n \,$ ($n \geq 2$).
In this case, $ W $ acts transitively on its roots,
so the algebras $ H_c $ and $  B_c $ are labelled
by a single parameter $ c \in \C\,$. We will write
$ \Cr $ for the set of regular values of $ c \,$, and will
use the standard notation denoting by $ \Q $ the field
of rational numbers and by $ \overline{\Q} $ its algebraic
closure in $ \C\,$. Note that $\, c \not\in  \Cr \,$ if and 
only if $\, e^{2 \pi i c} \,$ is a root of unity 
of order $\,d\,$, where $\, 2 \leq d  \leq n \,$ (see \cite{DJ}). 
Thus, we have $\, \Cr = \C \setminus \bigcup_{d=2}^{n} 
\bigcup_{m=1}^{d-1}\{m/d+ \Z\} \,$.

Our main results can be encapsulated into the following two theorems.
\begin{theorem} 
\label{TT1}
If $\, c \not\in \overline{\Q} \,$, the
algebras $ H_c $ and $ H_{c'} $ are

(a)\ isomorphic if and only if $ c = \pm c' \,$,

(b)\ $\C$-linearly Morita equivalent if and only if $ c \pm c' \in \Z \,$.

\end{theorem}
\begin{theorem} 
\label{TT2}
If $\,  c \not\in \overline{\Q} \,$, the
algebras $ B_c $ and $ B_{c'} $ are

(a)\ isomorphic if and only if $\, c = c' \,$ or $\, c = - c' - 1 \,$,

(b)\ $\C$-linearly Morita equivalent if and only if $ c \pm c' \in \Z \,$.
\end{theorem}

A few comments on the proof of these results.
First of all, the `if' parts of both statements in both theorems
are known. The isomorphism $\, H_c \cong H_{-c} \,$ is immediate
from the defining relations of $ H_c $ and holds for all 
$\, c \in \C\,$; $\, B_c \cong B_{-1-c} \,$ is a 
consequence of the fact that $ B_c $ can be identified with
the subalgebra of differential operators generated by $ \C[\h]^W $
and the Calogero-Moser operator $ L_c \,$, which is also true for 
all $\, c \in \C\,$ (see \cite{EG}, Proposition~4.9).
The existence of Morita equivalence 
between $\, H_c \,$ and $\, H_{c+k}\,$, $\, k \in \Z\,$, has been 
proven in \cite{BEG}, Theorem~8.1, for an arbitrary Coxeter group 
$ W\,$ and holds (at least) for regular $\,c\,$. A similar result 
for the spherical algebras is then immediate by Theorem~\ref{T0} above. 
Thus, we need only to establish the `only if' parts of 
Theorems~\ref{TT1} and~\ref{TT2}. We will prefer to work with
spherical algebras, proving first the required implications
in Theorem~\ref{TT2} (see Section~3 and Section~4 below). 
Theorem~\ref{TT1}$(b)$ will then follow immediately from  
Theorem~\ref{TT2}$(b)$ and Theorem~\ref{T0}, 
and with a little more work, we will conclude 
Theorem~\ref{TT1}$(a)$ from Theorem~\ref{TT1}$(b)$ (see Section~5).

The proof of Theorem~\ref{TT2} relies on explicit 
formulas for the Hatori-Stallings traces for the algebra $ B_c $
which we derive in Section~2. One general result that underlies
our computations is, perhaps, of independent interest and 
deserves to be mentioned here.
\begin{theorem} 
\label{TT3}
Let $ V $ be a finite-dimensional symplectic vector space over $ \C \,$, 
and let $\, G \subset \GL(V) \,$ be a finite group acting on $ V $
by symplectic transformations. Then the algebra  $\, A := \C[V]^G \,$ 
of polynomial $G$-invariants has a finite-dimensional $0$-th Poisson 
homology, i.e. $\, \dim_{\C\,}(A/ \{A, A\}) < \infty\,$,
where $\, \{A, A\} \subseteq A \,$ is the subspace spanned 
by Poisson brackets of all elements of $ A\,$.
\end{theorem}
\noindent The assertion of Theorem~\ref{TT3} has been conjectured
(and verified in many special cases) in \cite{AF}, \cite{AL}.
We will give a complete proof of this result in the Appendix.

We have stated Theorems~\ref{TT1} and~\ref{TT2} under the assumption
that one of the parameters ($ c $ or $ c' $) takes non-algebraic values. 
In fact, it suffices to prove these results only in the case when 
{\it both} $ c $ {\it and} $ c' $ are transcendental over $ \Q \,$. 
Once the latter is established, it is standard to conclude that the 
algebras $\, H_c \,$ (resp., $\, B_c $) and $\, H_{c'} \,$ 
(resp., $\, B_{c'} $) can be neither isomorphic nor Morita 
equivalent if one of the parameters is algebraic while the other 
is not. Indeed, if (say) $\, H_c \,$ were isomorphic to $ H_{c'} \,$
for some $\, c \in \overline{\Q} \,$ and $\, c'  \not\in \overline{\Q} \,$, 
then for {\it any} given $\, 
c'' \in \C \setminus \overline{\Q}\,$, we could find an automorphism
$\, g \in \Aut(\C\, / \,\overline{\Q})\,$, such that $\, g.c' = c'' \,$
and $ g.c = c \,$, forcing $ \, H_{c'} \cong H_{c''} \,$.
On the other hand, it seems very likely that both theorems remain 
true under milder assumptions. For example, it would be
interesting to extend the above classification to all
regular values of $ c \,$. Some partial results in this direction
are discussed in Section~6.

Finally, we should mention that in the simplest case ($ n=2 $), 
the algebras $\, B_c \,$ can be identified with primitive factors
of $U(\sll2)\,$ (see \cite{EG}, Section~8), and in this context 
both parts of Theorem~\ref{TT2} were proven earlier: $\, (a) $ is 
due to Dixmier \cite{D}, and $ (b) $ is due to Stafford \cite{S} 
and Hodges \cite{H}. The methods we employ in the present paper 
generalize those of Hodges.  In the end of \cite{H} 
the author remarks that $\K$-theoretic 
techniques can be applied to distinguish between primitive factors
of higher dimensional semisimple Lie algebras but ``... 
the question cannot be completely solved along this approach''. 
It seems surprising that in the case of Cherednik algebras, 
as we will see below, the $\K$-theoretic approach does lead 
to a complete solution of the problem, at least for generic 
parameter values.
\medskip

\noindent {\bf Acknowledgments.} {\footnotesize
We thank R.~Stanley for suggesting us the idea of the 
proof of Lemma~\ref{L2} in Section~4.
The first author was partially supported by the NSF grant 
DMS 00-71792 and A.~P.~Sloan Research Fellowship; the
work of the second author was  partly conducted for 
the Clay Mathematics Institute and partially
supported by the NSF grant DMS-9988796}.

\section{Hattori-Stallings Traces}
Let $ \Gamma = \Gamma(W) $ be the set of irreducible 
representations of $ W\,$.
The Grothendieck group of $ \C W $ is then
a free abelian group of rank $\, |\Gamma| \,$ generated by
the classes of $\, \tau \in \Gamma\,$:
$\, \K_0(\C W) = \bigoplus_{\tau \in \Gamma} \, \Z \cdot [\tau] 
\,$. According to \cite{EG}, the algebra $ H_c $ can be 
equipped with a canonical increasing filtration 
$ \{F_{\bullet} H_c \} \,$, whose first nonvanishing 
component $ F_0 H_c $ is $ \C W $ and the associated 
graded algebra $ \grd(H_c) $ is 
isomorphic to $ \C[\h \oplus \h^*] \# W \,$ (and therefore, is
Noetherian of finite global dimension). By a well-known theorem 
of Quillen \cite{Q}, the inclusion $ \C W \into H_c $ induces in 
this situation an isomorphism of abelian groups 
$ \K_0(\C W) \stackrel{\sim}{\to} \K_0(H_c) \,$. 
If we assume that $ c \in \Cr \,$, then 
$ B_c $ is Morita equivalent to $ H_c\,$ 
(see Theorem~\ref{T0} above), 
and the corresponding equivalence functor $\, 
\e H_c \otimes_{H_c} - : \Mod(H_c) \to \Mod(B_c) $ 
gives another isomorphism of $ \K$-groups: 
$\, \K_0(H_c) \stackrel{\sim}{\to} \K_0(B_c) \,$. 
Thus, for regular $ c \,$, we can identify both
$ \K_0(H_c) $ and $ \K_0(B_c) $ with $ \K_0(\C W)\,$.
Specifically, the class of an irreducible representation 
$ \tau \in \Gamma $ corresponds to the classes of (left) 
projective modules $ P_{\tau} := H_c \otimes_{\C W} \tau $
and $ \e P_{\tau} :=  \e H_c \otimes_{\C W} \tau $ in 
$ \K_0(H_c) $ and $ \K_0(B_c) $ respectively.
 
For any $ c\in \C\,$, $\, B_c $ is a Noetherian domain
(since so is $\, \grd(B_c) \cong \C[\h \oplus \h^*]^W \,$). 
By Goldie's theorem, it has a quotient 
division ring $ Q(B_c)\,$. Thus, we can define the rank 
homomorphism $\, \rk: \K_0(B_c) \to \Z  \,$ on $ \K_0(B_c) $ setting
\begin{equation}
\label{1}
\rk [M] := \mathtt{length}(Q(B_c) \otimes_{B_c} M)\ , \quad 
[M] \in \K_0(B_c) \ .
\end{equation}
Let $ \tilde{\K}_0(B_c) $ denote the kernel of (\ref{1}) 
so that $\, \K_0(B_c) \cong \tilde{\K}_0(B_c) \oplus \Z \,$.
Since $\, \rk[\e P_{\tau}] = \dim(\tau) \,$ for each 
$ \tau \in \Gamma\,$, we may take the classes
\begin{equation}
\label{2}
[\widetilde{\e P_{\tau}}] :=  \dim(\tau) \cdot [\e P_{\triv}] - 
[\e P_{\tau}] \ , \quad \tau \in \Gamma^{*}\ ,
\end{equation}
as generators of $ \tilde{\K}_0(B_c)\,$. 
(Here and below, $ \triv $ stands for the trivial representation of 
$ W $ and $ \Gamma^{*} :=   \Gamma \setminus \{\triv\}$.)

Now, let $ [H_c, H_c] $ be the subspace of $ H_c $ spanned 
by (additive) commutators, and let $ \HH_0(H_c) := H_c/[H_c, H_c] $ 
be the trace group (i.e., the $0$-th Hochschild homology) of 
the algebra $ H_c \,$ with canonical projection 
$\, \tr_{H_c}: H_c \onto  \HH_0(H_c)\,$.
Thus, we regard $ \HH_0(H_c) $ as a vector space over $ \C \,$, 
with $ \tr_{H_c} $ being a $\C$-linear map.

To study finite-dimensional representations of $ H_c $  we have 
calculated in \cite{BEG} the value of $ \tr_{H_c} $ 
on the central idempotents $ \e_{\tau} $ of  $ \C W \subset  H_c $ 
corresponding to simple modules $ \tau \in \Gamma\,$. 
The result of this calculation (see \cite{BEG}, Theorem~5.3) 
reads\footnote{\textbf{N.B.} Our notation here 
differs slightly from that of \cite{BEG}.}
\begin{equation}
\label{300} 
\tr_{H_c}(\e_{\tau}) = 
\frac{(\dim \,\tau)^2}{n! \, (nc)^n} F_{\tau}(nc)\cdot 
\tr_{H_c}(1)\ \quad \mbox{for any}\ c \not= 0 \ , 
\end{equation}
where $ F_{\tau}(x) \in \Z[x] $ is the so-called {\it content
polynomial} of the Young diagram $ Y(\tau) $ corresponding to 
$ \tau\,$ (see \cite{M}, Example~11, p. 15):
\begin{equation}
\label{310} 
F_{\tau}(x) := \prod\limits_{(i,j) \in Y(\tau)}(x + j-i) \ .
\end{equation}

Using (\ref{300}), it is easy to compute the {\it 
Hattori-Stallings trace} for the algebra $ H_c\,$:
$$
\ch_{H_c}:\, \K_0(H_c) \to \HH_0(H_c)\ .
$$
Indeed, by definition, if $ \e $ is an indempotent in $ H_c $, 
then $ H_c\e $ is a finite projective module whose Hattori-Stallings 
trace is given by $ \, \ch_{H_c}[H_c\e] := \tr_{H_c}(\e) \,$. Now, 
for any central group idempotent $\, \e_{\tau} \in \C W \,$, we have 
$\,[H_c\e_{\tau}] = \dim(\tau) \cdot [P_{\tau}] \,$ in $ \K_0(H_c) \,$.
Hence, it follows from (\ref{300}) that
\begin{equation}
\label{3}
\ch_{H_c}[P_\tau] = \frac{\dim(\tau)}{n! \,(nc)^n}\, 
F_{\tau}(nc)\cdot \tr_{H_c}(1)\ , \quad  c \not= 0\ .
\end{equation}

Now, write $\, \HH_0(B_c) := B_c/[B_c,B_c] \,$ for the trace group, 
$\, \tr_{B_c}: B_c \onto \HH_0(B_c) \, $ for the canonical projection, 
and $\, \ch_{B_c}:  \K_{0}(B_c) \to \HH_0(B_c) \,$ for  
the Hattori-Stallings trace of the spherical algebra $ B_c \,$. 
By definition, $\, \tr_{B_c}(1) $ is a value that 
$ \ch_{B_c} $ takes on the class $ [B_c] $ of the free 
module of rank one, i.e. $\, \ch_{B_c}[B_c] =  
\tr_{B_c}(1)\,$ in $\,\HH_0(B_c) \,$.
As we mentioned above, $ B_c $ is Morita equivalent to $ H_c \,$, 
when $ c $ is regular. Hence, in that case $ \HH_0(H_c) $ 
and  $ \HH_0(B_c) $ are isomorphic,
with isomorphism $\, \Phi_c: \HH_0(H_c) \stackrel{\sim}{\to}  
\HH_0(B_c) \,$ being naturally induced by the Morita functor. 
Moreover, the map $ \Phi_c $ fits in the commutative diagram
\begin{equation}
\label{410}
\begin{CD}
\K_{0}(H_c) @>\ch_{H_c}>> \HH_0(H_c)\\
@V \e VV     @VV\Phi_{c} V \\
\K_{0}(B_c) @>\ch_{B_c} >> \HH_0(B_c)\\
\end{CD}
\end{equation}
Since under the Morita equivalence $ [B_c] $ corresponds to the 
class of $ P_{\triv} $ in $ \K_{0}(H_c) $ and 
$\, \ch_{H_c}[P_{\triv}] = \frac{1}{n! \, (nc)^n} F_{\triv}(nc) 
\cdot \tr_{H_c}(1) \,$ by (\ref{3}), we see from (\ref{410}) that
\begin{equation}
\label{4}
\Phi_c(\tr_{H_c}(1)) = 
\frac{n! \, (nc)^n}{F_{\triv}(nc)}\cdot \tr_{B_c}(1) \ .
\end{equation}
Now, using (\ref{3}), (\ref{410}), (\ref{4}) and 
$ \C$-linearity of $ \Phi_c\,$, we compute
\begin{equation}
\label{5}
\ch_{B_c}[\e P_\tau] = \Phi_c(\ch_{H_c}[P_\tau]) =
\frac{\dim(\tau)\,F_{\tau}(nc)}{F_{\triv}(nc)}\cdot 
\tr_{B_c}(1)\ , \quad c \in \Cr \ .
\end{equation}
In particular, evaluating  $ \ch_{B_c} $ on generators of 
$ \tilde{\K}_{0}(B_c) \,$ (see (\ref{2})) yields
\begin{equation}
\label{6}
\ch_{B_c}[\widetilde{\e P_{\tau}}] = 
G_{\tau}(nc) \cdot \tr_{B_c}(1)\ , 
\quad \tau \in \Gamma^{*}\ , \ c \in \Cr \ ,
\end{equation}
where we denote
\begin{equation}
\label{600}
G_{\tau}(x) := \dim(\tau)\,\left(1 - \frac{F_{\tau}(x)}{F_{\triv}(x)}
\right) \ \in\  \mathbb{Q}(x)\ .
\end{equation}

The following result will be crucial for our further considerations.
\begin{proposition}
\label{P1}
We have

$ (a) $\ $\, \dim\,\HH_0(B_c) \geq 1 \,$ for all $ c \in \C \,$, 
and $\, \dim\,\HH_0(B_c) = 1 \,$ for all but finite subset of 
values of $ c\,$; 

$(b)$\ if $\, c \in \Cr \,$, then $\, \dim\, \HH_0(B_c) = 1 \,$,
and $\, \tr_{B_c}(1) \not= 0 \,$ in $\, \HH_0(B_c) \,$.
\end{proposition}
\noindent  We will deduce Proposition~\ref{P1} 
from Theorem~\ref{TT3} (which, in turn, will be proven in the Appendix). 
As a first step, we prove a weaker version of this result.
\begin{lemma}
\label{L0}
If $\, c \in \C \setminus \overline{\Q}\,$ then 
$\, \dim\,\HH_0(B_c) = 1 \,$, and 
$\, \tr_{B_c}(1) \not= 0 \,$ in $\, \HH_0(B_c) \,$.
\end{lemma}
\begin{proof}
Due to Morita invariance of 
Hochschild homology and formula (\ref{4}), 
we may replace the algebra $ B_c $ in the statement of the lemma 
by $ H_c \,$. Then, $\, \HH_{0}(H_c) \,$ being $1$-dimensional follows 
from the fact that, for generic $ c \,$ $\, \HH_{0}(H_c) \,$ is 
isomorphic to the 0-th Poisson homology of the Calogero-Moser 
space $ \mathcal{M}_n\,$ which, in turn, is isomorphic to the 
top (De Rham) cohomology of $ \mathcal{M}_n\,$ (see \cite{EG}). 
By a theorem of Nakajima (\cite{N}, Theorem~3.45), $\,\mathcal{M}_n\,$ is 
diffeomorphic (as a $ \mathcal{C}^{\infty}$-manifold) 
to the Hilbert scheme $ \Hilb_{n}(\C^2) $ of $n$-points on $ \C^2 \,$, the top Betti 
number of which is known to be equal to one (see \cite{N}, 
Corollary~5.10).

A proof of the fact that $\, \tr_{H_c}(1) \not= 0 \,$ for generic $ c\,$ 
was sketched in \cite{BEG} (cf. Remark after Theorem~5.3, 
{\it loc. cit.}). 
For reader's convenience, we recall and slightly refine 
this argument here. 
Let $\, E := \{\,c \in \C\, :\, 1 \in [H_c, H_c] \ \mbox{in}\ H_c\,\}\,$.
Then, for each $\, c \in E \,$, there exists 
$\, i = i(c) \in \Z_{\geq 0} \,$, such that
$\,1 \in [F_i H_c, F_i H_c ] \,$, where 
$\, F_{\bullet}H_c \,$ is the standard increasing 
filtration on $ H_c\,$. Thus, we have $\, 
E = \bigcup_{i=0}^{\infty} E_i \,$, where 
$\, E_i := \{c \in \C\, :\, 1 \in [F_i H_c, F_i H_c]\}\,$
for $\, i = 0, 1, 2 \ldots $ Now, each $ E_i $ is
a semi-algebraic set defined over $ \Q\,$. On the other hand,
$\, E_i \,$ cannot contain a non-empty Zariski open
subset of $ \C\,$. Indeed, if this were the case for some $ i\,$,
then $\, \C \setminus E_i \,$ would be at most a finite set.
The latter is impossible, because, as we know (see \cite{BEG1}, 
Theorem~1.2), there are infinitely many values of $ c $ of 
 the form $\, c = 1/n + \Z_{\geq 0} \,$, for which the algebra
$ H_c $ has nonzero finite-dimensional representations
(the element $ 1 \in H_c $ acts as the identity operator on such 
a representation, and hence has a nonzero trace). 
Thus, each $\, E_i\,$ is a (possibly empty) algebraic subset of 
$ \C $ defined over $ \Q\,$, i.e.
$\, E_i \subset \overline{\Q}\,$ for all $\, i = 0, 1, 2 \ldots $ 
It follows that $\, E \subset \overline{\Q}\,$, 
and hence $\, \tr_{H_c}(1) \not= 0 \,$ when
$\, c \in \C \setminus \overline{\Q}\,$. This 
finishes the proof of the lemma.
\end{proof}

Now, to prove Proposition~\ref{P1} in full generality 
we introduce the following notation. First, we define 
an associative algebra $ \H $ to be generated by the 
elements of $\, \h \,$, 
$\, \h^* \,$, $\, W \,$, and a new {\it central variable} $\, 
\cc \,$ satisfying the same relations as $ H_c\,$. Next, we set 
$\, \B := \e \H \e \,$, where $ \e $ is the $W$-symmetrizer in $ \H\,$. 
By construction, the algebra $ \B $ (as well as $ \H $) has 
a non-trivial center, namely $\,\C[\cc]\,$, each $ B_c $ being a 
quotient of $ \B $ obtained by specializing a central character. 
Moreover, the trace group $\,\HH_{0}(\B) := \B/[\B,\B]\,$ has the
structure of a module over $\,\C[\cc]\,$, and we have
\begin{lemma}
\label{L8}
$ \HH_{0}(\B) $ is a finite  $\C[\cc]$-module.
\end{lemma}
\begin{proof}
Recall that each  $ B_c $ is equipped 
with the increasing filtration $\, \{F_{\bullet}B_c\} \,$, 
the associated graded ring
$\,\grd(B_c) \cong \C[\h \oplus \h^*]^W\,$ 
being independent of $ c\,$. Letting $ \deg(\cc)= 0\,$,
we can extend this filtration to the algebra $ \B \,$
so that $\, \grd(\B) \cong A[\cc]\,$. Then we have
$\,\grd(\B)/\{\grd(\B),\, \grd(\B)\} \cong (A/\{A,A\})[\cc]\,$,
where $ A:= \C[\h \oplus \h^*]^W \,$. 
By Theorem~\ref{TT3}, $\, A/\{A,A\} \,$ is a finite-dimensional
vector space over $ \C \,$. Hence, 
$\, \grd(\B)/\{\grd(\B),\, \grd(\B)\} \,$ is a finite module 
over $ \C[\cc]\,$. Now, if we equip $\, [\B,\B] \subseteq \B \,$ and
$\, \HH_{0}(\B) =  \B/[\B,\B] \,$ with the induced filtrations, 
then $\,\grd([\B,\B]) \subseteq  \grd(\B) \,$ and
$\,\grd\,\HH_{0}(\B) \cong \grd(\B)/\grd([\B,\B]) \,$. As
$\, \{\grd(\B),\,\grd(\B)\} \subseteq \grd([\B,\B]) \,$, we see
that $\,\grd\,\HH_{0}(\B)\,$ is a quotient of 
$\, \grd(\B)/\{\grd(\B),\,\grd(\B)\} \,$, and hence, is finite over 
$\, \C[\cc]\,$. This implies that 
$\, \HH_{0}(\B) \,$ is finite over $ \C[\cc]\,$ as well. 
Lemma~\ref{L8} is proven.
\end{proof}

Now we are in position to give a complete

\vspace*{1ex}

\noindent
{\sl Proof of Proposition~\ref{P1}}. In view of Lemma~\ref{L8}, 
we may think of $ \HH_{0}(\B) $ as a {\it coherent} 
sheaf on $ \Spec\,\C[\cc]\,$, with $\, \HH_{0}(B_c) \,$ being a 
fiber over the point $\, c \,$. The function $\, c \mapsto 
\dim\,\HH_{0}(B_c) \,$ is then upper semicontinuous 
(see, e.g., \cite{Ha}, Exercise~II.5.8) which means 
that, for all $ m \in \Z\,$, the sets $\, \{c \in \C\,:\, 
\dim\,\HH_{0}(B_c) \geq m \}\,$ are (Zariski) closed in 
$\C\,$. With Lemma~\ref{L0}, this gives
immediately the first statement of our proposition: 
to wit, $\, \dim\,\HH_{0}(B_c) \geq 1 $ for all $ c \,$,
the set $\, \{c \in \C\,:\, \dim\,\HH_{0}(B_c) > 1 \}\,$ 
being at most finite. Now, if $\, c \in \Cr \,$, then 
$\, \HH_{0}(B_{c+k}) \cong \HH_{0}(B_c) \,$ for any $\, k \in \Z \,$ 
(since $ B_{c+k} $ is Morita equivalent to $ B_c \,$; 
see \cite{BEG}, Theorem~8.1), and therefore we have necessarily
$\, \dim\,\HH_{0}(B_c) = 1 \,$ in that case. Next, by Lemma~\ref{L0}, 
$\, \tr_{B_c}(1) \not= 0 \,$ for transcendental values of $ c\,$.
Being a section of the coherent sheaf $ \HH_{0}(\B) \,$, 
$\,\tr_{B_c}(1)\,$ may vanish then only on a closed, and therefore {\it finite} 
subset in $ \C\,$. Again, the vanishing of $\, \tr_{B_c}(1) $ for a regular 
value of $\, c \,$ would imply that $\, \tr_{B_{c+k}}(1) = 0\,$ for all $\, k \in \Z \,$ 
(cf. (\ref{16}) and (\ref{17}) below). Thus, to avoid a contradiction
we conclude that $\, \tr_{B_c}(1) \not= 0 \,$ when $\, c \in \Cr \,$.
Proposition~\ref{P1} is proven. 

\vspace{2ex}

\begin{remark}
\label{RR1} Note that the second statement of 
Proposition~\ref{P1}$(b)$ does not hold for the algebra $ H_c $ 
(even though Lemma~\ref{L0} does). Indeed, if $ c=0 $ then 
$\, H_c \cong \dd(\h) \# W \,$ and $ \tr_{H_c}(1) = 0 \,$. We expect, 
however, that $ c=0 $ is the only exceptional value, and 
$ \tr_{H_c}(1) $ does not vanish for any $ c \not= 0\,$. On the other
hand, we also expect that $\, \dim\, \HH_{0}(H_c) = 1 \,$ for all 
values of $c\,$. In the simplest case, when $ n=2 $ and $ B_c $
are identified with primitive quotients of $U(\sll2)\,$, 
$ \HH_{0}(B_c) $ (and, more generally, $ \HH_{0}(B_{c}^G) $ for all
finite $ G \subset \Aut(B_c) $) have been computed explicitly in \cite{F}.
\end{remark}

\section{Proof of Theorem~\ref{TT2}$(a)$}
The Hattori-Stallings traces have good functorial properties 
with respect to change of rings and Morita equivalence
(see, e.g, \cite{B}). 
Specifically, given an algebra homomorphism 
$ \varphi: B \to B' \,$, there is a commutative diagram 
\begin{equation}
\label{7}
\begin{CD}
\K_{0}(B) @>{\ch_{B}}>> \HH_0(B)\\
@V {\K_{0}(\varphi)}VV     @VV {\HH_0(\varphi)} V \\
\K_{0}(B') @>{\ch_{B'}} >> \HH_0(B')\\
\end{CD}
\end{equation}
where $ \K_{0}(\varphi) $ is a homomorphism of abelian groups 
sending $ [P] \in \K_{0}(B) $ to the class of the induced 
module $ [B' \otimes_{B} P] \in \K_{0}(B') \,$ and
$\, \HH_0(\varphi): \HH_0(B) \to \HH_0(B')$ is a $ \C$-linear map 
given by $\, b + [B,B]\, \mapsto \, \varphi(b) + [B',B'] \,$
(see \cite{B}, Section~2).

We will use (\ref{7}) to distinguish two generic members 
of the family of algebras $ \{B_c\} $ up to isomoprhism. 
First, we observe that if $ \varphi: B_c \to B_{c'} $ is
an isomorphism of $ \C$-algebras, then $ \varphi(1) = 1 \,$, 
and therefore, by definition of $ \HH_0(\varphi) $, we have
\begin{equation}
\label{700}
\HH_0(\varphi)(\tr_{B_c}(1)) = \tr_{B_{c'}}(1)\ .
\end{equation}
Next, the map $ \K_{0}(\varphi) $ preserves rank, and hence
restricts to an isomorphism of reduced $\K$-groups:
$ \tilde{\K}_{0}(\varphi): \tilde{\K}_{0}(B_c) \to 
\tilde{\K}_{0}(B_{c'}) \,$. Choosing the classes 
$\, [\widetilde{\e P_{\tau}}(c)] \,$ and 
$\, [\widetilde{\e P_{\tau}}(c')] \,$ (see (\ref{2}))
as bases in $ \tilde{\K}_{0}(B_c) $ and $ \tilde{\K}_{0}(B_{c'}) $
respectively, we can represent the isomorphism 
$ \tilde{\K}_{0}(\varphi) $ by an invertible integral-valued 
matrix $\, \| m_{\tau \sigma} \| 
\in \GL_{N}(\Z) \,$ of dimension $ N = |\Gamma| - 1\,$:
$$
\tilde{\K}_{0}(\varphi)\,:\ [\widetilde{\e P_{\tau}}(c)] 
\, \mapsto \, \sum_{\sigma \in \Gamma^*} 
m_{\tau \sigma}\cdot [\widetilde{\e P_{\sigma}}(c')]
\ , \quad \tau \in \Gamma^* \ .
$$
The commutative diagram 
\begin{equation}
\label{9}
\begin{CD}
\tilde{\K}_{0}(B_c) @>{\ch_{B_c}}>> \HH_0(B_c)\\
@V {\tilde{\K}_{0}(\varphi)}VV     @VV{\HH_0(\varphi)}V  \\
\tilde{\K}_{0}(B_{c'}) @>{\ch_{B_{c'}}} >> \HH_0(B_{c'})\\
\end{CD}
\end{equation}
produces a system of equations
$$
\sum_{\sigma \in \Gamma^*} m_{\tau \sigma}\cdot  
\ch_{B_{c'}}[\widetilde{\e P_{\tau}}(c')] = 
\left(\HH_0(\varphi) \,\ccirc\, \ch_{B_c}\right)
[\widetilde{\e P_{\tau}}(c)] \ , 
\quad \tau \in \Gamma^* \ ,
$$
which can be written explicitly 
(use (\ref{6}), (\ref{600}) and (\ref{700})) as follows
\begin{equation}
\label{10}
\left(\sum_{\sigma \in \Gamma^*} 
m_{\tau \sigma}\,G_{\sigma}(nc') - G_{\tau}(nc) \right) \cdot \tr_{B_{c'}}(1) = 0\ .
\end{equation}

Now we are in position to prove the first part of 
Theorem~\ref{TT2}. 
As mentioned in the Introduction, we need only to establish 
the implication: $\, B_c \cong B_{c'}\, \Rightarrow \, c=c'\ 
\mbox{or}\ c =-c'-1 \,$, and for $ n = 2 \,$, this is already
known (due to Dixmier \cite{D}, Th\'eor\`eme~6.4, and
Hodges \cite{H}, Theorem~3).

Thus, we fix $\, n \geq 3 \,$ and assume
that $\, B_{c} \cong B_{c'} \,$ for some $\, c \not \in \overline{\Q}\,$.
Then, as explained in the Introduction, we may also assume
that $\, c' \not \in \overline{\Q}\,$. 
By Lemma~\ref{L0}, we have $\, \tr_{B_{c'}}(1) \not= 0 \,$, and
hence, for $\, c, c' \not\in \overline{\Q}\,$, the equations 
(\ref{10}) are satisfied if and only if
\begin{equation}
\label{111}
\sum_{\sigma \in \Gamma^{*}} m_{\tau \sigma}\,
G_{\sigma}(nc') = G_{\tau}(nc)\ ,\quad \tau \in \Gamma^{*} \ .
\end{equation}

Now, the trivial representation of $ W $ corresponds to the 
partition $ (n) $ so that $\, F_{\triv}(x) = \prod_{k=0}^{n-1}(x + k)\,$ 
by (\ref{310}). Each function $\, G_{\lambda}(x) \,$, $\, \lambda \in
\Gamma^* \,$, can be developed then into elementary fractions 

\begin{equation}
\label{12}
G_{\lambda}(x) = \sum\limits_{k=1}^{n-1} \frac{a_{\lambda, k}}{x+k}
\end{equation}
with $\, a_{\lambda,k} \in \Q \,$. Specifically,
we see from (\ref{600}) that
\begin{equation}
\label{13}
a_{\lambda,k} = - \dim(\lambda)\, \frac{\prod\limits_{(i,j) 
\in Y(\lambda)}(j-i-k)}
{\prod\limits_{{l=0 \atop l\not=k}}^{n-1}(l - k)}\ .
\end{equation}
Taking the subset of representations 
$ \{\lambda_m \, |\, m =1, 2, \ldots,  n-1 \} $ indexed 
by the partitions $ \Par(\lambda_m) := (m, 1^{n-m})\,$, 
we check easily (with (\ref{13})) that $ a_{\lambda_m, k} = 0 $ 
for $ k < m $, while $ a_{\lambda_m, m} \not= 0 $ for all 
$ m = 1,2, \ldots, n-1 \,$. Hence, each of elementary 
fractions $\, 1/(x+k)\,$ that occur in 
(\ref{12}) can be expressed as a $ \mathbb{Q}$-linear combination 
of the functions $\, G_{\lambda_1}(x), G_{\lambda_2}(x), \ldots, 
G_{\lambda_{n-1}}(x) \,$. 
It follows then from (\ref{111}) that there are some numbers
$\, b_{kj} \in \mathbb{Q} \,$ (depending on $ m_{\tau \sigma} $) 
such that
$$
\frac{1}{nc+k} = \sum_{j=1}^{n-1}\frac{b_{kj}}{nc'+j}\quad 
\mbox{for}\ k = 1, 2, \ldots, n-1 \ .
$$
Letting $ x = nc' \,$, we can rewrite these equations in the form
\begin{equation}
\label{14}
nc = \frac{f(x)}{g_1(x)} - 1 = 
\frac{f(x)}{g_2(x)} - 2 = \ldots =  
\frac{f(x)}{g_{n-1}(x)} - (n-1)\ ,
\end{equation}
where $\, f(x) := (x+1)(x+2)\ldots (x+n-1) \,$ and 
$\, g_1(x),\, g_2(x), \, \ldots,\, g_{n-1}(x) \,$ are some nonzero 
polynomials in $ \mathbb{Q}[x] $ of degree $ \leq  n-2 \,$.

Under the assumption that $ c' \not\in \overline{\Q}\,$, 
all the equalities in (\ref{14}) starting with the second one, 
should hold as identities 
in $ \mathbb{Q}(x)\,$ (that is, not only for $ x = nc' \,$ 
but for {\it all} $\, x \in \C\,$). 
Indeed, if this were not the case for some $ k \,$, the corresponding 
difference $\, f(x)/g_{k+1}(x) - f(x)/g_{k}(x) - 1 \,$ would 
provide a non-trivial rational polynomial having $ c' $ as a root. 
Now, since $ \deg[f(x)] > \deg[g_k(x)] \,$, the zero set of each 
function $ f(x)/g_{k}(x) $ is a non-empty subset of 
$ \{-1,-2, \ldots, 1-n\} \,$, the set of roots of $\, f(x) \,$. 
On the other hand, any two of these functions, 
say $\, f(x)/g_{j}(x) $ and  $ f(x)/g_{k}(x) \,$ 
with $ j \not= k \,$, cannot have zeros in common because 
$ f(x)/g_{j}(x) - f(x)/g_{k}(x) = j - k  $ by (\ref{14}).
Since the number of fractions $ f(x)/g_{k}(x) $ matches
exactly the number of zeros of $ f(x)\,$, we conclude that 
each $ f(x)/g_{k}(x) $ has {\it one and only one} zero, 
and therefore, must be of the form
$\, f(x)/g_{k}(x) = \frac{1}{q_k} (x - x_k) \,$ with some 
$ q_k \in \Q^{\times} \,$ and $\, 
x_k \in \{-1,-2, \ldots, 1-n\} \,$. Substituting these into (\ref{14}) 
gives $\, q_1 =  q_2 = \ldots = q_{n-1} =: q \,$ and 
$\, x_k = x_1 + (1-k)q \,$, where $ k = 1, 2, \ldots, n-1 \,$.
Since each $ x_k \in \{-1,-2, \ldots, 1-n\} \,$, we may have
only two possibilities: $\, x_k = -k,\, q = 1 \,$ and 
$\, x_k = k-n, \,q = -1 \,$, which give the two required 
relations: $\,  c = c' $ and  $  c = -c'-1 \,$ respectively.
This finishes the proof of Theorem~\ref{TT2}$(a)$.

\section{Proof of Theorem~\ref{TT2}$(b)$}
Our next goal is to classify the algebras $ \{ B_c \}$ up to 
Morita equivalence. For this, we will extend and slightly refine 
the argument given in the proof of Theorem~\ref{TT2}$(a)$. We start 
with a general (and perhaps, well-known) ring-theoretic 
result to be needed later. 
\begin{lemma}
\label{L1}
Let $ B $ be a Noetherian domain, and let $ P $ be a 
finitely-generated projective module which is 
a generator in the category of right $B$-modules. 
Assume that $ B' := \End_{B}(P) $ is a domain. Then,
$ P_{B} $ is isomorphic to a right ideal in $ B\,$, 
and $ {}_{B'}P $ is isomorphic to a left ideal in $ B'\,$ .
\end{lemma}
\begin{proof} By Morita's theorem, $ B' $ is equivalent 
to $ B\,$, and hence is Noetherian. By Goldie's theorem, 
$ B' $ satisfies then  Ore's condition: every pair
of nonzero right ideals in $ B' $ has nonzero intersection. 
This means that $ B' \,$ (regarded as
a right module over itself) is uniform. Being uniform
is a Morita-invariant property (see \cite{MR}, Lemma~3.5.8(vi)). 
Hence, the $ B$-module $ P $ which corresponds to 
$ B' $ under the Morita equivalence is also uniform.

Now, by the Dual basis lemma (see, e.g., \cite{MR}, Lemma~3.5.2), 
we have $ P \cdot P^* = B'\,$, where $ P^* := \Hom_{B}(P, B)\,$.
Therefore we may find $\, \xi_1, \ldots, \xi_n \in P \,$
and $ \theta_1, \ldots, \theta_n \in \ms{P}^*\,$,
such that $ \sum_{i=1}^{n}\xi_i \cdot \theta_i = 1 \,$. 
We claim that at least one of the maps $ \theta_i $ is injective. 
To see this, first observe that $ \bigcap_{i=1}^{n}\,\Ker(\theta_i) = 
\{0\} \,$. Indeed, if $ \xi \in \bigcap_{i=1}^{n}\,
\Ker(\theta_i) \,$, then 
$\xi = 1.\xi = \sum_{i=1}^{n} \xi_{i} \theta_{i}(\xi) = 0\,$.
Now, assuming that $ \Ker(\theta_i) \not= \{0\} $ for 
each $ i = 1, 2, \ldots, n \,$, we may choose $ m < n $
so that $ \bigcap_{i=1}^{m}\,\Ker(\theta_i) \not= \{0\} $
while $ \bigcap_{i=1}^{m+1}\,\Ker(\theta_i) = \{0\} \,$.
Hence,
$\,\bigcap_{i=1}^{m}\,\Ker(\theta_i) \,\bigoplus\, 
\Ker(\theta_{m+1}) \subseteq P\,$,
which contradicts the uniformity of $ P\,$.
Thus, $ \theta_i\,: P \to B $ is injective for some
$ i \,$, and we may use it to identify $ P $ with a right 
ideal in  $ B \,$. 

Finally, by the standard Morita theory (cf. \cite{MR}, 
Corollary 3.5.4(b)), if $ P_B $ is a finite projective
generator in the category of right $B$-modules,
then $ {}_{B'}P $ is a finite projective
generator in the category of left $B'$-modules, and 
$ \End_{B'}({}_{B'}P) \cong B\,$. Repeating the above 
argument verbatim (with roles of $B$ and $B'$ interchanged)
shows that $ {}_{B'}P $ is isomorphic to a left ideal in $B'\,$.
\end{proof}

Now, we recall that any equivalence functor $ \ms{F}:\, \Mod(B) \to
\Mod(B') \,$ between module categories is isomorphic to 
$ P \otimes_B \mbox{---} \,$ for some finitely generated 
projective module $ P = P_{B}\,$.
Moreover, such an $ \ms{F} $ induces isomoprhisms
$ \K_0(\ms{F}): \K_0(B) \to \K_0(B') $ and 
$ \HH_0(\ms{F}): \HH_0(B) \to \HH_0(B') $ making commutative
the following diagram
\begin{equation}
\label{15}
\begin{CD}
\K_{0}(B) @>{\ch_{B}}>> \HH_0(B)\\
@V {\K_{0}(\ms{F})}VV     @VV {\HH_0(\ms{F})} V \\
\K_{0}(B') @>{\ch_{B'}} >> \HH_0(B')\\
\end{CD}
\end{equation}
By definition, the map $ \K_0(\ms{F}) $ takes the class of 
the free module $ [B] \in \K_0(B) $ to $ [{}_{B'}P] \in \K_0(B') \,$, 
while $ \ch_{B}[B] = \tr_{B}(1)\,$ in $ \HH_0(B)\,$. Hence, by 
commutativity of (\ref{15}), we have
\begin{equation}
\label{16}
\HH_0(\ms{F})(\tr_{B}(1)) = \ch_{B'}[P]\ .
\end{equation}
Returning to our situation, let $ B := B_c \, , 
\, B' := B_{c'} \,$, and let $ \ms{F} $
be a $ \C$-linear equivalence of categories: $\,\Mod(B_{c}) \to
\Mod(B_{c'})\,$. Since both $ B_c $ and $ B_{c'} $ are domains,
we have $\, \ms{F} \simeq  P \otimes_{B_c} \mbox{---} \,$,
with  $ P = {}_{B'}P_{B} $ being isomorphic to a right projective 
ideal in $ B_c $ and to a left projective ideal in $ B_{c'} $
(see Lemma~\ref{L1}). Regarding $ P $ as the latter,
we have $\, \rk[P] = \rk[B_{c'}] = 1\,$, and therefore $\, [P] - [B_{c'}] \in 
\tilde{\K}_{0}(B_{c'})\,$. We can write then 
$\,[P] = [B_{c'}] + \sum_{\lambda \in \Gamma^*} n_{\lambda} 
\cdot [\widetilde{\e P_{\lambda}}(c')]\,$ in $ \K_{0}(B_{c'}) $
for some $ n_{\lambda} \in \Z\,$, and compute (with (\ref{6}))
\begin{equation}
\label{17}
\ch_{B_{c'}}[P] = \biggl(1 + \sum_{\lambda \in \Gamma^*} 
n_{\lambda}\, 
G_{\lambda}(nc')\biggr) \cdot \tr_{B_{c'}}(1)\ , \quad c' \in \Cr \ .
\end{equation}
Looking at the diagram (\ref{15}) and taking into account 
(\ref{16}) and (\ref{17}), we find
\begin{eqnarray}
(\HH_0(\ms{F}) \,\ccirc \, \ch_{B_{c}})[\widetilde{\e P_{\tau}}(c)]
&=& \HH_0(\ms{F})(G_{\tau}(nc) \cdot \tr_{B_{c}}(1))
\nonumber\\*[2ex]
&=&  G_{\tau}(nc)\, \HH_0(\ms{F})(\tr_{B_{c}}(1)) \qquad (\mbox{by
$\C$-linearity of $ \ms{F} $})\nonumber\\
&=& G_{\tau}(nc)\,\biggl(1 + \sum_{\lambda \in \Gamma^*} n_{\lambda}\, 
G_{\lambda}(nc')\biggr) \cdot \tr_{B_{c'}}(1) \nonumber
\end{eqnarray}
On the other hand,
\begin{eqnarray}
(\ch_{B_{c'}}\, \ccirc \, \K_0(\ms{F}))[\widetilde{\e P_{\tau}}(c)]
&=& \sum_{\sigma \in \Gamma^*} m_{\tau \sigma} \cdot 
\ch_{B_{c'}}[\widetilde{\e P_{\sigma}}(c')]\nonumber\\
&=&  \sum_{\sigma \in \Gamma^*} m_{\tau \sigma}\,
G_{\sigma}(nc') \cdot \tr_{B_{c'}}(1) \ ,\nonumber
\end{eqnarray}
where $ \|  m_{\tau \sigma}\| \in \GL_{N}(\Z) $ is an invertible 
integral-valued matrix representing the isomorphism $ \K_0(\ms{F}) $ 
in the bases 
$\,\{[\widetilde{\e P_{\sigma}}(c)]\} \subset \K_0(B_c) \,$ 
and $\,\{[\widetilde{\e P_{\sigma}}(c')]\} \subset \K_0(B_{c'}) \,$.
Thus, by commutativity of (\ref{15}), we have 
\begin{equation}
\label{18}
\left(\sum_{\sigma \in \Gamma^*} m_{\tau \sigma}\,
G_{\sigma}(nc') - G_{\tau}(nc)\,
\biggl(1 + \sum_{\lambda \in \Gamma^*} n_{\lambda}\, 
G_{\lambda}(nc')\biggr)\right)\cdot \tr_{B_{c'}}(1) = 0\ .
\end{equation}
When $\, c' \not\in \overline{\Q}\,$, we have
$\, \tr_{B_{c'}}(1) \not= 0 \,$,
so in that case (\ref{18}) is equivalent to the system of equations
\begin{equation}
\label{19}
\sum_{\sigma \in \Gamma^*} m_{\tau \sigma}\,
G_{\sigma}(nc') = 
G_{\tau}(nc)\,\biggl(1 + \sum_{\lambda \in \Gamma^*} n_{\lambda}\, 
G_{\lambda}(nc')\biggr)\ , \quad  \tau \in \Gamma^*\ .
\end{equation}

Now, we are in position to prove part $ (b) $  of Theorem~\ref{TT2}.
First of all, we observe that the same argument as we used 
in the proof of the first part of this theorem (reducing 
the system (\ref{111}) to (\ref{14})) works for the system 
(\ref{19}) as well. 
As a result, we get from (\ref{19}) a set of $ (n-1) $ algebraic
equations of the form (\ref{14}),
with $ f(x) $ depending now on 
$  a_{\lambda, k} $ and new parameters $ n_{\lambda} 
\in \Z \,$. Specifically, 
\begin{equation}
\label{20}
f(x) = \prod\limits_{k=1}^{n-1} (x+k) + 
\sum\limits_{k=1}^{n-1}\, a_{k}\,
\prod\limits_{{j=1 \atop j\not=k}}^{n-1}(x+j)\ ,
\end{equation}
where $\, a_{k} := \sum_{\lambda \in \Gamma^*} n_{\lambda} 
\, a_{\lambda, k} \,$ (and $  a_{\lambda, k} $ are defined by
(\ref{12}), (\ref{13})). 
Under the assumption that $ c' $ is non-algebraic, the system
(\ref{14}) again gives us the relations
\begin{equation}
\label{21}
nc = (x-x_1)/q-1 = (x-x_2)/q-2 = \ldots = (x-x_{n-1})/q-(n-1)\ ,
\end{equation}
where $ \{x_1, x_2, \ldots, x_{n-1}\} $ is the set of roots
of (\ref{20}) and $\, q \in \Q^{\times}\,$. Summing up these 
relations and using Vi\`ete's theorem for $ f(x) \,$, 
we find 
$$ 
(n-1)nc = (n-1)nc'/q + 
\sum_{k=1}^{n-1}(k+a_{k})/q - n(n-1)/2\ ,
$$
or equivalently,
\begin{equation}
\label{22}
q\biggl(c+ \frac{1}{2}\biggr) = \biggl(c' + \frac{1}{2}\biggr) + 
\frac{1}{n(n-1)}\sum\limits_{k=1}^{n-1}\, a_{k} \ .
\end{equation}
Thus, to finish our proof it remains
to show that $\, q = \pm 1 \,$ and the last term in
(\ref{22}) is an integer. As we will see, the former follows
from the latter, while the latter is an immediate 
consequence of the following observation.
\begin{lemma}
\label{L2}
For any $\, \lambda \in \Gamma^* \,$ and $\, k = 1, 2, 
\ldots, n-1 \,$, the numbers $ a_{\lambda, k} $ are integers 
divisible by $ n(n-1)\,$.
\end{lemma}
Indeed, by (\ref{13}), we have
\begin{equation}
\label{23}
a_{\lambda, k} = (-1)^{n-k-1} \,
\frac{\dim (\lambda)}{k!\ (n-1-k)!}\, F_{\lambda'}(k) \ ,\quad 
k = 1, 2, \ldots, n-1 \ ,
\end{equation}
where $ F_{\lambda'}(x) $ is the content polynomial of the Young
diagram conjugate to $ Y(\lambda)\,$. Using the ``hook formula'' 
for the dimension of irreducible representations of $ S_n $ 
(see, e.g., \cite{J}, Theorem~20.1) and
a well-known formula for the Schur function (\cite{M}, Example~4,
p.~45), we can rewrite (\ref{23}) in the form 
\begin{equation}
\label{24}
a_{\lambda, k} = (-1)^{n-k-1} \, n \, {n-1 \choose k}\, 
\sch_{\lambda'}(1^k)\ ,
\end{equation}
where $\, \sch_{\lambda'}(1^k) \,$ is the Schur function 
of the conjugate partition of $ \lambda $ evaluated at 
$ \boldsymbol{\mathtt{x}} = (1, \ldots,1 , 0,\ldots,0) \,$ 
(the first $ k $ symmetric variables are equal to $ 1 \,$, 
the rest are zero). So we need only to see that 
$\, {n-1 \choose k}\, \sch_{\lambda'}(1^k) \,$ is divisible
by $\, n-1\,$. Since each Schur function
$ \sch_{\lambda'}(\boldsymbol{\mathtt{x}}) $ can be 
written as a linear combination of monomial 
symmetric functions $ \mon_{\sigma}(\boldsymbol{\mathtt{x}}) $
with integer coefficients\footnote{These 
coefficients are usually called {\it Kostka numbers} 
(see \cite{M}, Chapter~I, Section~6).}:
\begin{equation}
\label{240}
\sch_{\lambda'}(\boldsymbol{\mathtt{x}}) = \sum
\limits_{\sigma \in \Gamma} K_{\lambda' \sigma}\,
\mon_{\sigma}(\boldsymbol{\mathtt{x}})\ ,\quad  
K_{\lambda' \sigma} \in \Z\ ,
\end{equation}
it suffices to show that $\, 
\frac{1}{n-1}\,{n-1 \choose k}\, \mon_{\sigma}(1^k) \,$
is an integer for every $ \sigma \in \Gamma\,$. By 
\cite{M}, Example~1(a), p.~26, we have
\begin{equation}
\label{25}
{n-1 \choose k}\, \mon_{\sigma}(1^k) = 
{n-1 \choose k}\ \frac{l!}
{\mu_{1}!\ \mu_{2}!\ \ldots }\ {k \choose l}
= \frac{(n-1)!} {(n-1-k)!\ (k-l)!\
\mu_{1}!\ \mu_{2}!\ \ldots } \ ,
\end{equation}
where $ l = l(\sigma) $ is the length of the partition of
$ \sigma $ and $ \mu_{i} = \mu_{i}(\sigma) $ is the number 
of parts of $ \Par(\sigma) $ equal to $ i \,$. Now, 
since $\, (n-1-k)+ (k-l) + \mu_1 + \mu_2 + \ldots = n-1 \,$
while $ \sum_{i\geq 1} i\,\mu_{i} = n \,$, the numbers
$\,\{n-1-k,\, k-l,\, \mu_1,\,  \mu_2,\, \ldots \} \,$ are relatively 
prime. Hence, there exist some $\, a, b, c, d,\ldots \in \Z \,$ 
such that 
$$
(n-1-k)\,a + (k-l)\,b + \mu_1\, c +  \mu_2 \,d + \ldots = 1 \ .
$$
Multiplying both sides of this equation by the multinomial
coefficient of the right-hand side of (\ref{25}), we see 
that the latter is divisible by $ n-1\,$. This finishes
the proof of Lemma~\ref{L2}.
%

Next, we must show that $ q = \pm 1 $ in (\ref{22}).
By (\ref{21}), we have $\, x_k  = x_1 + (1-k)q\,$ for 
$\, k = 1, 2, \ldots, n-1\,$, and since $\, q \in \Q \,$ and 
$\, \sum x_k \in \Z \,$, all the roots of $ f(x) $ are 
rational.  But, in view of Lemma~\ref{L2}, $ f(x) $ is a monic polynomial 
with {\it integer} coefficients. Hence, being rational, 
all $ x_k$'s are, in fact, in $ \Z\,$. 
It follows that $ q = x_1 - x_2  \in \Z\,$. Now, in order to 
conclude that $ q $ equals $ \pm 1 \,$, 
we simply reverse the roles of $ c $ and $ c' $ in the above 
considerations. More precisely, we have seen already 
that $\, q\,(c+1/2) - (c'+1/2) \in \Z \,$
with $ q \in \Z \,$; hence, by symmetry, 
we must also have $\, q'\,(c'+1/2) - (c+1/2) \in \Z \,$ with 
some $ q' \in \Z\,$. It follows that $\, (qq'-1)\,(c'+1/2) \in \Z \,$, 
and since we assume $ c' $ to be transcendental over $ \Q\,$, this 
implies $\, q\,q' =1 \,$, i.e. $ q $ is a unit in $ \Z\,$.
With Lemma~\ref{L2}, this argument completes the 
proof of the theorem.

\vspace*{2ex}

\begin{remark}
\label{R1}
Using formulas (\ref{24}), (\ref{240}) and (\ref{25}), we can
compute the sum in (\ref{22}) explicitly:
\begin{eqnarray} 
\sum_{k=1}^{n-1}\, a_{k} 
&=& \sum_{k=1}^{n-1} \ \sum_{\lambda \in \Gamma^*}\, n_{\lambda}\,a_{\lambda, k} 
= \sum_{k=1}^{n-1}\, (-1)^{n-k-1}\, n \,{n-1 \choose k} \
\sum_{\lambda \in \Gamma^*}\, n_{\lambda}\,\sch_{\lambda'}(1^k)\nonumber\\ 
&=&
\sum_{k=1}^{n-1} \,(-1)^{n-k-1}\, n \,{n-1 \choose k}\ 
\sum_{\lambda \in \Gamma^*}\, n_{\lambda}\ 
\sum_{\sigma \in \Gamma} \, K_{\lambda' \sigma} \, \mon_{\sigma}(1^k) \nonumber\\
&=&
\sum_{\sigma \in \Gamma}\, \tilde{n}_{\sigma} \ 
\sum_{k=1}^{n-1} \,(-1)^{n-k-1} \, n \,{n-1 \choose k}\,
\mon_{\sigma}(1^k)\nonumber\\
&=&
\sum_{\sigma \in \Gamma}\, \tilde{n}_{\sigma}\ 
\sum_{k=1}^{n-1} \,(-1)^{n-k-1} \,\frac{n!}{(n-1-k)!\ (k-l)!\ 
\mu_{1}!\ \mu_{2}!\ \ldots}\nonumber\\ 
&=&
\sum_{\sigma \in \Gamma}\, \tilde{n}_{\sigma}\, 
\frac{n\,(n-1)\,\ldots \,(n-l)}{\mu_{1}!\ \mu_{2}!\ \ldots}\
\sum_{k=1}^{n-1}\, (-1)^{n-k-1}\, \frac{(n-l-1)!}{(n-1-k)!\ (k-l)!}\ ,
\nonumber
\end{eqnarray}
where we denote $\, \tilde{n}_{\sigma} := 
\sum_{\lambda \in \Gamma^*}\, K_{\lambda' \sigma}\, n_{\lambda} \,$
for each $ \sigma \in \Gamma \,$. By Newton's binomial formula, 
the last sum is zero unless $\, l = n-1 \,$, and it is equal to 
$ 1 $ in that case. Now there is only one partition of weight 
$ n $ and length $ l = n-1 \,$, namely $\, (2, 1^{n-2}) \,$, 
in which case we have $\, \mu_{1} = n-2 \, ,\, \mu_2 = 1\,$ and 
$ \mu_i = 0 \,$ for $\, i \geq 2\,$. Thus,
\begin{equation}
\label{26}
\sum_{k=1}^{n-1}\, a_{k} =\, n\,(n-1)\, \tilde{n}_{\alpha}\, =\,
n\,(n-1)\,\sum_{\lambda \in \Gamma^*}\, K_{\lambda' \alpha}\ 
n_{\lambda}\ ,
\end{equation}
where $\, \alpha \in \Gamma \,$ is the representation 
of $ W $ corresponding to the partition $ (2, 1^{n-2})\,$.
Now, comparing (\ref{22}) and (\ref{26}), we see that 
the parameters of Morita equivalent algebras $ B_c $ and 
$ B_{c'} $ are related to the trace of the corresponding 
equivalence functor (\ref{17}) by the formula
$$
\pm(c+1/2) = (c'+1/2) +
\sum_{\lambda \in \Gamma^*}\, K_{\lambda' \alpha}\ n_{\lambda}\ .
$$
\end{remark}

\section{Proof of Theorem~\ref{TT1}}

In view of Theorem~\ref{T0}, 
Morita classifications for the families of algebras $ \{H_c\} $ 
and $ \{B_c\} $ 
must be equivalent, at least in case when $ c \in \Cr\,$. 
Thus, part $(b)$ of Theorem~\ref{TT1} follows from part 
$(b)$ of Theorem~\ref{TT2} which we have proved in 
Section~4 above.
In part $(a)$, the `if' implication
is obvious, and we only need to show its converse, i.e.
$\, H_c \cong H_{c'} \, \Rightarrow \, c = \pm c' \,$.

First of all, when isomorphic, the algebras $ H_c $ and $ H_{c'} $ 
are Morita equivalent, and therefore, by part $(b)\,$, we have at once
\begin{equation}
\label{27}
c = \pm c' + l \quad \mbox{for some}\ l \in \Z\ .
\end{equation}
Arguing now as in Section~3 and using
the trace formula (\ref{3}) for $ H_c \,$,
we may derive a system of equations similar to (\ref{10}),
namely
\begin{equation}
\label{28}
\left(\frac{\dim(\tau)}{n! \,(nc)^n}\,F_{\tau}(nc) -
\sum_{\sigma \in \Gamma} \, m_{\tau \sigma}\, 
\frac{\dim(\sigma)}{n! \,(nc')^n}\, 
F_{\sigma}(nc')\right)\cdot \tr_{H_{c'}}(1) = 0 \ , \quad \tau \in \Gamma\ ,
\end{equation}
where $\, \| m_{\tau \sigma}\| \in \GL_{|\Gamma|}(\Z) \,$
and $\, F_{\sigma}\,,\, F_{\tau} \,$ are defined in (\ref{310}).
Again, as $\, c' \not\in \overline{\Q} \,$, the trace 
factor $\, \tr_{H_{c'}}(1) \,$ can be dropped by Lemma~\ref{L0},
and we can rewrite (\ref{28}) in the form
$$
F_{\tau}(nc)\,(nc')^n = 
\sum_{\sigma \in \Gamma}\, 
k_{\tau \sigma}\, F_{\sigma}(nc')\,(nc)^n \ ,\quad k_{\tau \sigma} \in 
\Q \ .
$$
In particular, letting $\, \tau = \triv \,$ and taking into account
(\ref{27}), we have
\begin{equation}
\label{29}
(\pm x + nl + 1)\ldots(\pm x + nl + n-1)\,x^n = 
\sum_{\sigma \in \Gamma}\, 
k_{\tau \sigma}\, F_{\sigma}(x)\,(\pm x + nl)^{n-1} \ ,
\end{equation}
where $ x = nc'\,$. Since $ c' \not\in \overline{\Q}\,$,
the equation (\ref{29}) should hold identically in $ x\,$. 
Substituting then $\, x = \mp nl \,$ gives 
$\,(n-1)!\,(\mp nl)^{n} = 0 \,$, whence $\, l = 0\,$.
With (\ref{27}), this finishes the proof of 
Theorem~\ref{TT1}$(a)$.

\section{Concluding Remarks}

As already mentioned in the Introduction, we expect 
Theorems~\ref{TT1} and~\ref{TT2} to be true not only for 
transcendental but for all regular values of $ c\,$. 
This is immediate in case $ n=2 $ (see \cite{H}) but 
seems to be much harder to prove in higher dimensions.
Reviewing the proofs of Section~3 and~4 shows that their 
argument works also for algebraic $ c \,$, provided the 
degree of the extension of fields $\,\Q(c)/\Q\,$ is large 
enough: specifically, $\, [\Q(c) : \Q] \geq  n(n-1)\,$.

It is worth noting that the field $ \Q(c) $ is a Morita
invariant for the whole family of algebras $ \{B_c\} \,$ 
(with no restrictions imposed on $c$).
Indeed, if both $ c $ and $ c' $ are regular, then using
Proposition~\ref{P1}$(b)$ (instead of Lemma~\ref{L0})
still allows one to reduce $ (\ref{18}) $ to the form 
(\ref{14}), and the first equation of (\ref{14}) implies 
immediately that $\, \Q(c) \subseteq \Q(c')\,$. 
By symmetry, we then also have $\, \Q(c') \subseteq \Q(c)\,$.
Now, if $\, c \,,\, c'\,$ are both singular then 
$\, c \,,\, c' \in \Q \,$, and therefore $\,\Q(c)=\Q(c')=\Q\,$. 
On the other hand, if one of the parameters is singular, say $\, c \,$, 
while $\, c' \,$ is regular, the algebras $ B_c $ and $ B_{c'} $ 
cannot be Morita equivalent. Indeed, in this case 
$ B_{c'} $ is a simple ring (by Theorem~\ref{T0}) while
$ B_{c} $ is not. To see the latter, note that if $\, c \not\in 
\Cr \,$ then $ c $ is singular in the sense of \cite{DJO}\,:
\,for such $ c \,$, the standard $H_c$-module $\, M(\triv) = 
\C[\h] \,$ corresponding to the trivial representation is reducible.
Let $ J  $ be a proper submodule of $\, M(\triv) \,$ which we regard 
as an ideal in $ \C[\h] \,$. Then, $\, J^W := J\, \cap \, \C[\h]^W \,$ 
is a nonzero ideal in $ \C[\h]^W \,$ acting trivially on the nonzero $ B_c$-module 
$ \e(M/J)\,$. Thus, $\, J^W \,$ is contained in the annihilator
of $ \e(M/J) \,$, and hence $\, \Ann_{B_c}[\e(M/J)] \,$ 
is a proper two-sided ideal of $ B_c \,$.

In the end, we would like to make one, perhaps somewhat speculative 
observation. According to Theorem~\ref{TT1}, a Morita class 
of $\, H_c \,$ is a function (at least, for $ c \not\in \overline{\Q} $) 
of $\, e^{2 \pi ic}\,$, that is, of the  monodromy representation 
of the system of Dunkl operators associated to $\, (W,c)\,$ (see \cite{Du}). 

\vspace*{1ex}

\question
Is there a genuine relation between the Morita classification of 
Cherednik algebras $ H_c(W) $ and the monodromy representation of 
the corresponding Hecke algebras $ \mathcal{H}_{W}(e^{2 \pi ic})\,$?

\section{Appendix: Proof of Theorem~\ref{TT3}}

The purpose of this Appendix is to prove Theorem~\ref{TT3} 
announced in the Introduction. Our proof below does not rely
on and can be read independently of the rest of the paper.

First, we fix notation and recall the claim to be proved. 
Let $ V $ be a complex symplectic vector space, i.e.
a finite-dimensional vector space over $ \C $ equipped
with a nondegenerate skew-symmetric bilinear form 
$ (\,\boldsymbol{\cdot} \,,\, \boldsymbol{\cdot}\,)\,$.
Let $\, G \subset \SSp(V) \,$ be a finite group of 
linear symplectic transformations of $ V \,$. 
Write $ A := \C[V]^{G} $ for the algebra of polynomial 
$ G$-invariants. The form 
$ (\,\boldsymbol{\cdot} \,,\, \boldsymbol{\cdot}\,)\,$
determines a structure 
of the (graded) Poisson algebra on $ \C[V] \,$, and this 
structure restricts naturally to $ A \,$. Let $\,\{A,A\} \,$ denote
the linear span of Poisson brackets of all elements of $ A\,$. 
Then Theorem~\ref{TT3} states: $\,\dim_{\C\,}(A/\{A,A\}) < \infty\,$.

To show this, observe that the quotient $ A/\{A,A\} $ 
has a natural grading, each graded component 
$\, (A/\{A,A\})_n := A_n/\sum_{i+j=n+2}\{A_i, A_j\} \,$ 
being finite-dimensional. Set $\, 
(A/\{A,A\})^{*}_{n} := \Hom_{\C\,}\left[\,(A/\{A,A\})_n\,,\,\C\,\right]\,$
for $ n = 0, 1, 2,\ldots\,$ Then we have
\begin{lemma}
\label{LA1}
$ (A/\{A,A\})^{*}_{n} $ is isomorphic 
to the subspace of polynomials $\, P \in \C[V] \,$ of degree 
$ n $ satisfying the equations
\begin{equation}
\label{E1}
\sum_{g \in G}\, (u, gv)\, P(u + gv) = 0 \quad 
\mbox{for all}\  u, v \in V\ .
\end{equation}
\end{lemma}
\begin{proof}
For a fixed $\, u \in V \,$, define the function 
$\, L_u\,:\, V \to \C\,$, $\, x \mapsto L_u(x) := \sum_{g \in G}\, 
e^{(u,gx)} \,$. We can think of $\, L_{u}(x) \,$ as a 
generating function for the graded vector space $ A\,$ 
(meaning that the coefficients of $ L_{u}(x) \,$, 
when $ L_{u}(x) $ is regarded as a series in $u\,$, span $ A\,$). 
Now, observe that $\,\{ L_u, L_v\}(x) = 
\sum_{g \in G}\, (u,gv)\, L_{u+gv}(x) \,$ for all 
$\, u, v \in V\,$. Hence, if $\,f: A \to \C \,$ is a linear 
homogeneous functional such that $\, f(\{A,A\})=0\,$, then
$\, u \mapsto f(L_u)\,$ is a homogeneous $G$-invariant
polynomial $ P(u) \,$ satisfying (\ref{E1}).
Conversely, if $ P $ satisfies (\ref{E1}) then $ P $ defines a 
functional $\, f_{P} \,$ on $ A $ vanishing on $ \{A,A\}\,$.
\end{proof}

Thus, to prove Theorem~\ref{TT3}, it suffices to show that 
the functional equations (\ref{E1}) may have only finitely many 
linearly independent polynomial solutions. Replacing 
$ v \mapsto v t $ and expanding the left-hand side of (\ref{E1}) 
into Taylor series in $ t \,$, we get a system of linear 
differential equations for $ P\,$:
\begin{equation}
\label{E2}
\sum_{g \in G}\, (u, gv)\, \partial_{gv}^{m} P(u) = 0 \quad 
\mbox{for all}\ v \in V\ \mbox{and}\ m = 0, 1, 2, \ldots\ ,
\end{equation}
where $\, \partial_{gv}^{m} P(u) := (d/dt)^m P(u + gvt)|_{t=0} \,$.
It suffices to prove that the space of local (holomorphic) 
solutions of (\ref{E2}) in a neighborhood of some point 
$\, u_0 \in V\, $ is finite-dimensional. This would, in turn, 
follow if we establish that $\,\C[V]\,$ is a finite module over 
the subalgebra the symbols of differential equations (\ref{E2}) 
generate at $ u_0\,$. Thus, we need to show that $\,\C[V]\,$ is 
finite over the algebra generated by polynomials 
$\, x \mapsto Q_{v,m}(x) := \sum_{g \in G}\, (u_0, gv)\,(x,gv)^m
\,$ for all $\, v \,$ and $\, m \,$.

Assume the contrary. Then, being homogeneous, the ideal of $ \C[V] $ 
generated by $\, Q_{v,m}\,$ must be of infinite codimension in 
$ \C[V] \,$, and therefore vanishes at some point $ x_0 \not= 0\,$.
In other words, the system of polynomial equations 
$\, Q_{v,m}(x) = 0 \,$, and therefore the system
\begin{equation}
\label{E3}
\sum_{g \in G}\, (u_0, gv)\, e^{t(x,gv)} = 0 \quad 
(\mbox{for all}\ v \in V \ \mbox{and}\ t \in \C)
\end{equation}
have a nonzero solution $ x = x_0\,$.

Let $\, O_{x_0} \subset V \,$ be the orbit, and 
$\, G_{x_0} \subseteq G \,$ be the stabilizer of $ x_0 $ under the 
action of $ G\,$. Choose $ v \in V $ in such a way that 
$\, (x,v) \not= (x',v) \,$ whenever $\, x, x' \in O_{x_0} \,$ 
and $\,  x \not= x' \,$. It follows then from (\ref{E3}) that 
$\,(u_0, \sum_{g \in G_{x_0}}\,gv) = \sum_{g \in G_{x_0}}\, 
(u_0, gv) = 0 \,$. In other words, $ u_0 $ is orthogonal to 
$\, \sum_{g \in G_{x_0}}\,gv \,$ for generic (and hence, 
any) vector $ v \in V\,$. Thus, $ u_0 $ is orthogonal to the 
image of the linear $ G_{x_0}$-symmetrizer acting on $ V \,$.
Every such image is a nonzero subspace of $ V $ (as it contains 
$ x_0 $), and there are only finitely many of these. So, if a priori 
we choose $\, u_0 \,$ not to be orthogonal to any of such subspaces, 
we get a contradiction finishing the proof of the theorem.

\vspace*{1.5ex}

We close this section with two remarks related to the work 
of J.~Alev et al. (see \cite{AF}, \cite{AL}).

\vspace*{1ex}

\remark It has been shown in \cite{AF} that 
$\, \dim_{\C\,}(A/\{A,A\}^2) < \infty \,$, where $ \{A,A\}^2 $ 
is the subspace of $ A $ spanned by all elements of the form 
$\,\{a,b\}\{c,d\}\,$. This result follows from Theorem~\ref{TT3}\,:
more generally, the theorem implies that $\, \{A,A\}^n \,$ has 
finite codimension in $ A $ for every $ n \geq 2\,$. 
Indeed, being graded and of finite codimension, the subspace 
$\,\{A,A\}\,$ contains a power of the augmentation ideal 
of $ A\,$, i.e. $\, (A_{+})^{m} \subseteq \{A,A\} \,$ for 
some $ m \geq 1\,$. Hence, $\, (A_{+})^{mn} \subseteq \{A,A\}^n \,$,
and therefore $\, \{A,A\}^n \,$ has also finite codimension in 
$ A\,$.

\vspace*{1ex}

\remark  In Section~2, we apply Theorem~\ref{TT3} in the 
situation when $\, V = \h \oplus \h^* \,$ and $ G = W $ 
is the symmetric group (so that $\, A = \C[\h \oplus \h^*]^W \,$). 
In this particular case we expect that $\, \dim_{\C\,}(A/\{A,A\}) = 1 \,$. 
This would follow from Proposition~\ref{P1} (and Th\'eorem\`e~6 of 
\cite{AL}), if one would be able to show that $\,\HH_0(B_c)\,$ is a 
{\it deformation} of the Poisson homology of 
$\, \grd(B_c)\,$, i.e.
$$
\grd\,\HH_0(B_c) = \grd(B_c)/\{\grd(B_c)\,,\,\grd(B_c)\}\ ,
$$
where $\,\grd\,\HH_0(B_c) \,$ is formed with respect to the quotient 
induced filtration on $\,\HH_0(B_c)\,$. An alternative way would 
be to solve exactly the system of functional equations (\ref{E1}) 
in case when $\, V = \h \oplus \h^* \,$ and $ G = W \,$.

\footnotesize{

\end{document}